\documentclass[12pt]{article}

\usepackage[utf8]{inputenc} 
\usepackage[T1]{fontenc}    
\usepackage{hyperref}       
\usepackage{url}            
\usepackage{booktabs}       
\usepackage{amsfonts}       
\usepackage{nicefrac}       
\usepackage{microtype}      

\usepackage{graphicx}
\usepackage[square, comma, numbers, sort&compress]{natbib}
\usepackage{doi}
\usepackage{esdiff}

\usepackage[utf8]{inputenc}
\usepackage[english]{babel}
\usepackage[utf8]{inputenc}
\usepackage[T1]{fontenc}
\usepackage{lmodern}
\usepackage{color}

\usepackage{amsthm}
\usepackage{graphicx}
\usepackage{float}
\usepackage{mathtools,amssymb,dsfont,IEEEtrantools}
\usepackage{amsmath}

\usepackage[verbose=true,letterpaper]{geometry}
\AtBeginDocument{
  \newgeometry{
    textheight=9in,
    textwidth=6.5in,
    top=1in,
    headheight=14pt,
    headsep=25pt,
    footskip=30pt
  }
}

\usepackage{enumitem}
\setlist{listparindent=\parindent,parsep=0pt}

\newtheorem{theorem}{Theorem}[section]
\newtheorem{lemma}[theorem]{Lemma}

\theoremstyle{definition}
\newtheorem{definition}[theorem]{Definition}

\theoremstyle{remark}
\newtheorem{remark}[theorem]{Remark}

\hypersetup{
pdftitle={A small remark on Bernstein's theorem},
pdfsubject={},
pdfauthor={Michael Bildhauer, Bernhard Farquhar, Martin Fuchs},
pdfkeywords={},
}

\title{A small remark on Bernstein's theorem}
\author{Michael Bildhauer\\ \url{bibi@math.uni-sb.de} \and Bernhard Farquhar\\ \url{farquhar@math.uni-sb.de} \and Martin Fuchs\\ \url{fuchs@math.uni-sb.de}}
\date{ Saarland University\\ 
        Department of Mathematics\\
        66123 Saarbrücken, Germany}

\renewcommand{\newline}{}
    
\newcommand{\dx}{\, \mathrm{d}x}

\begin{document}
\maketitle

\begin{center}
    \today
\end{center}

\begin{abstract}
\noindent We investigate splitting-type variational problems with some linear growth conditions. For balanced solutions of the associated 
Euler-Lagrange equation we receive a result analogous to Bernstein's theorem on non-parametric minimal surfaces.

Without assumptions of this type,  Bernstein's theorem cannot be carried over to the splitting case, which follows from an elementary 
counterexample. 

We also include some modifications of our main theorem.%
\footnote{AMS Classification 49Q20, 49Q05, 53A10, 35J20}, 
\end{abstract}

\section{Introduction}
	A famous theorem of Bernstein (see \cite{MR1544873}) states that a smooth solution $u=u(x), x=(x_{1},x_{2})$, 
	of the non-parametric minimal surface equation 
	\begin{equation}
        \mathrm{div}\left(\frac{\nabla u}{\sqrt{1+|\nabla u|^2}}\right)=0 \label{1}
    \end{equation}
    defined on the whole plane must be an affine function. Letting $f_{0}(P)\vcentcolon=\sqrt{1+|P|^2},P\in\mathbb{R}^2$, 
    the validity of \eqref{1} on some domain $\Omega\subset\mathbb{R}^2$ just expresses the fact that $u$ is a solution of the Euler-Lagrange 
    equation
	 associated to the area functional
    \begin{equation}
        J_{0}[u,\Omega]\vcentcolon=\int_{\Omega}f_{0}(\nabla u) \, \mathrm{d}x \label{2}.
    \end{equation}
    For a general overview on minimal surfaces, variational integrals with linear growth and for a careful analysis of 
    Bernstein's theorem the reader is referred for instance to  \cite{MR2778928}, \cite{MR2760441}, \cite{MR2566897}, 
    \cite{MR1015936}, \cite{MR90833} and  \cite{MR852409} and the references quoted therein.\\
    
    \newline 
    We ask the following question: does Bernstein's theorem extend to the case when the area integrand $f_{0}(\nabla u)$ is 
    replaced by the energy density $\sqrt{1+(\partial_{1}u)^2}+\sqrt{1+(\partial_{2}u)^2}$ being also of linear growth with 
    respect to $|\nabla u|$ but without any obvious geometric meaning? \\
    
    \newline 
    More generally we let for $P=(p_{1},p_{2})\in\mathbb{R}^2$
    \begin{equation}
        f(P):=f_{1}(p_{1})+f_{2}(p_{2}) \label{3}
    \end{equation}
    with functions $f_{i}\in \mathit{C}^2(\mathbb{R}), i=1,2$, satisfying
    \begin{equation}
        0< f_{i}''(t)\leq C_{i}(1+t^2)^{-\frac{\mu_{i}}{2}}, t\in\mathbb{R}, \label{4}
    \end{equation}
    for numbers $C_{i}>0$ and with exponents
    \begin{equation}
        \mu_{i}>1. \label{5}
    \end{equation}
    
    Note that \eqref{4} implies the strict convexity of $f$ and on account of \eqref{5} the density $f$ is of linear growth in the sense that
    \begin{equation*}
        |f_i(t)|\leq a|t|+b, \quad t\in\mathbb{R},
    \end{equation*}
    for some constants $a,b>0$. For a discussion of the properties of densities $f$ satisfying \eqref{3}-\eqref{5} we refer to \cite{MR4151290}.
     We then replace \eqref{1} by the equation
    \begin{equation}
        \mathrm{div}\left(Df(\nabla u)\right)=0 \label{6}
    \end{equation}
    and observe that the non-affine function ($\alpha,\beta,\gamma,\delta \in\mathbb{R}$)
    \begin{equation}
        w(x_1,x_2):=\alpha x_1x_2+\beta x_1 +\gamma x_2 +\delta \label{21}
    \end{equation}
    is an entire solution of equation \eqref{6}, in other words: the classical version of Bernstein's theorem does not extend 
    to the splitting case. The behaviour of the function $w$ defined in \eqref{21} is characterized in
    
    \begin{definition}
    A function $u\in\mathit{C}^{1}(\mathbb{R}^2), u=u(x_{1},x_{2})$, is called unbalanced, if and only if both of the following conditions hold
                  \begin{eqnarray}
                \limsup_{|x|\to \infty}\frac{|\partial_{1}u(x)|}{1+|\partial_{2}u(x)|}&=&\infty, \label{7}\\
                 \limsup_{|x|\to \infty}\frac{|\partial_{2}u(x)|}{1+|\partial_{1}u(x)|}&=& \infty. \label{8}
            \end{eqnarray}
           
        Otherwise we say that $u$ is of balanced form.
    \end{definition}
    
    \begin{remark}
    Condition \eqref{7} for example means that there exists a sequence of points $x_{n}\in\mathbb{R}^2$ such that 
    $|x_{n}|\to\infty$ and for which $$\lim_{n\to\infty}\frac{|\partial_{1}u(x_{n})|}{1+|\partial_{2}u(x_{n})|}=\infty .$$
    \end{remark}
    
    \begin{remark}
    If for instance  \eqref{7} is violated, no such sequence exists. Thus we can find constants $R,M>0$ such that 
    $|\partial_{1}u(x)|\leq M\left(1+|\partial_{2}u(x)|\right)$ for all $|x|\geq R$. Since $u$ is of class 
    $\mathit{C}^1(\mathbb{R}^2)$, this just shows
    
    \begin{equation}
        |\partial_{1}u(x)|\leq m|\partial_{2}u(x)|+M,\quad x\in\mathbb{R}^2, \label{9}
    \end{equation}
    with suitable new constants $m,M>0$.
    
    \end{remark}
Now we can state the appropriate version of Bernstein's theorem in the above setting:
\begin{theorem}
Let \eqref{3}-\eqref{5} hold and let $u\in\mathit{C}^2(\mathbb{R}^2)$ denote a solution of \eqref{6} on the entire plane. 
Then $u$ is an affine function or of unbalanced type.
\end{theorem}

\begin{remark}
    We do not know if \eqref{21} is the only entire unbalanced solution of \eqref{6}.\\
\end{remark}

Before proving Theorem 1.4 we formulate some related results:
in Theorem 1.7 below we can slightly improve the result of Theorem 1.4 by adjusting the notation introduced in \mbox{Definition 1.1} 
and by taking care of the growth rates of the second derivatives $f''_{i}$ (compare \eqref{4}).

\begin{definition}
    Let $\mu:=(\mu_1,\mu_2)$ with numbers $\mu_{i}>1$, $i=1,2$. A function $u\in\mathit{C}^1(\mathbb{R}^2)$ 
    is called $\mu$-balanced, if we can find a positive constant $c$ and a number $\rho>0$ such that at least one of the 
    following inequalities holds:
    
            \begin{equation}
                |\partial_1 u|\leq c(|\partial_2 u|^{\rho\mu_2}+1), \label{24}
            \end{equation}
        
            \begin{equation}
                |\partial_2 u|\leq c(|\partial_1 u|^{\rho\mu_1}+1), \label{25}
            \end{equation}

    where in case \eqref{24} we require $\rho\in(1/\mu_2,1)$, whereas in case \eqref{25} $\rho\in(1/\mu_1,1)$ must hold.
\end{definition}
Note that for example \eqref{24} is a weaker condition in comparison to \eqref{9}.\\

The extension of Theorem 1.4 reads as follows:
\begin{theorem}
    Let \eqref{3}-\eqref{5} hold and let $u\in\mathit{C}^2(\mathbb{R}^2)$ denote an entire solution of \eqref{6}. 
    If the function $u$ is $\mu$-balanced, then it must be affine.
\end{theorem}

In Theorem 1.8 we suppose that $|\partial_1u|$ is controlled in $x_2$-direction and from this we derive a 
smallness condition in $x_1$-direction -- at least for a suitable sequence satisfying $|x_1|\to\infty$. 
The idea of proving Theorem 1.8 again is of Bernstein-type in the sense that the proof follows the ideas of Theorem 1.4 
combined with a splitting structure of the test functions.\\

\begin{theorem}
    Let \eqref{3}-\eqref{5} hold and let $u\in C^2(\mathbb{R}^2)$ denote a solution of \eqref{6} on the entire plane. \newline
    Suppose that there exist real numbers $\kappa_1>0$, $ 0\leq\kappa_2<1$ such that
    \begin{equation}
        \mu_1>1+\frac{1}{\kappa_1}\frac{\kappa_2}{1-\kappa_2} \label{27}
    \end{equation} 
    and such that with a constant $k>0$
    \begin{equation}
        \limsup_{|x_2|\to\infty}\frac{|\partial_1u(x_1,x_2)|}{|x_2|^{\kappa_2}}\leq k. \label{28}
    \end{equation}
    Then we have
    \begin{equation}
        \liminf_{|x_1|\to\infty}\frac{|\partial_1 u(x_1,x_2)|}{|x_1|^{\kappa_1}}=0. \label{29}
    \end{equation}
    More precisely, by \eqref{27} we choose $\rho$ such that
    \begin{equation}
        \frac{\kappa_2}{1-\kappa_2}<\rho<\kappa_1(\mu_1-1). \label{30}
    \end{equation}
    Then the $\limsup$ in \eqref{28} is taken in the set
    \begin{equation*}
        M_{2,\rho}\vcentcolon=\left\lbrace(x_1,x_2) \colon |x_1|\leq2|x_2|^{1/(1+\rho)}\right\rbrace
    \end{equation*}
    and the $\liminf$ in \eqref{29} is taken in the set
    \begin{equation*}
        M_{1,\rho}\vcentcolon=\left\lbrace(x_1,x_2) \colon |x_2|\leq2|x_1|^{1+\rho}\right\rbrace.
    \end{equation*}
\end{theorem}

Our final Bernstein-type result is given in Theorem 1.9. Here a formulation in terms of the densities $f_i$ is presented without
requiring an upper bound for the second derivatives $f_i''$ in terms of some negative powers (see \eqref{4}).

\begin{theorem}
    Suppose that $f_i \in \mathit{C}^2(\mathbb{R})$, $i=1,2$, satisfies $f_i''(t) > 0$ for all $t\in \mathbb{R}$ and
    $f_i' \in L^{\infty}(\mathbb{R})$. Let $u\in\mathit{C}^2(\mathbb{R}^2)$ denote an entire solution of \eqref{6}, i.e.~it holds
    \[
    0= \mathrm{div}\left(Df(\nabla u)\right) = f_1''(\partial_1 u) \partial_{11}u + f_2''(\partial_2 u)\partial_{22}u \quad\mbox{on $\mathbb{R}^2$}.
    \]
    If 
    \[
    \Theta := \frac{f_2''(\partial_2 u)}{f_1''(\partial_1 u)} \in L^\infty(\mathbb{R}^2) \quad
    \mbox{or}\quad \frac{1}{\Theta} \in L^\infty (\mathbb{R}^2),
    \]
    then $u$ is an affine function.\\
\end{theorem}

In the next section we prove our main Theorem 1.4 while in Section 3 the variants mentioned above are established.

\section{Proof of Theorem 1.4}

Our arguments make essential use of a Caccioppoli-type inequality involving negative exponents. 
This result was first introduced in \cite{MR4151290}.We refer to the presentation given in Section 6 of \cite{MR4502898}, 
where Proposition 6.1 applies to the situation at hand. Let us assume that the conditions \eqref{3}-\eqref{5} hold 
and that $u$ is an entire solution of equation \eqref{6} being not necessarily of balanced type.

\begin{lemma}[see \cite{MR4502898}, Prop.6.1]
Fix $l\in\mathbb{N}$ and suppose that $\eta\in\mathit{C}_{0}^{\infty}(\Omega)$, $0\leq\eta\leq 1$, where $\Omega$ 
is a domain in $\mathbb{R}^2$. Then the inequality 
\begin{eqnarray}
    \lefteqn{\int_{\Omega}D^{2}f(\nabla u)(\nabla\partial_{i}u,\nabla\partial_{i}u)\eta^{2l}\Gamma_{i}^{\alpha}\, \mathrm{d}x} \nonumber\\
    &\leq &\int_{\Omega}D^{2}f(\nabla u)(\nabla\eta,\nabla\eta)\eta^{2l-2}\Gamma_{i}^{\alpha+1}\,
    \mathrm{d}x, \quad \Gamma_{i}\vcentcolon=1+|\partial_{i}u|^2 ,\label{22}
\end{eqnarray}
holds for any $\alpha >-1/2$ and for any fixed $i=1,2$.
\end{lemma}

Here and in what follows the letter $c$ denotes finite positive constants whose value may vary from line to 
line but being independent of the radius. \\

\newline
Assume next that the solution $u$ is balanced and w.l.o.g. let $u$ satisfy \eqref{9}. 
In order to show that $u$ is affine we return to inequality \eqref{22}, choose $i=1$ and fix some function 
$\eta\in C^{\infty}_0(B_{2R}(0))$ according to $\eta\equiv 1$ on $B_{R}(0)$, $|\nabla\eta|\leq c/R$.
Then \eqref{22} yields for any exponent $\alpha\in(-1/2,\infty)$ and with the choice $l=1$ ($B_r := B_r(0)$, $r>0$)
\begin{eqnarray}
 \lefteqn{\int_{B_{2R}}D^2f(\nabla u)(\nabla\partial_{1}u,\nabla\partial_{1}u)\eta^2\Gamma^{\alpha}_{1}\, \mathrm{d}x}\nonumber\\
 &\leq&  c\int_{B_{2R}}D^2f(\nabla u)(\nabla\eta,\nabla\eta)\Gamma_{1}^{\alpha+1}\,\mathrm{d}x \nonumber \\
 &\overset{\eqref{3}}{=}&c\int_{B_{2R}-B_{R}}\left(f_{1}''(\partial_{1}u)|\partial_1\eta|^2
 +f_{2}''(\partial_{2}u)|\partial_2\eta|^2\right)\Gamma_{1}^{\alpha+1}\,\mathrm{d}x\nonumber \\
 &\leq& cR^{-2}\left(\int_{B_{2R}-B_{R}}f_{1}''(\partial_{1}u)\Gamma_{1}^{\alpha+1}\,\mathrm{d}x
 +\int_{B_{2R}-B_{R}}f_{2}''(\partial_{2}u)\Gamma_{1}^{\alpha+1}\,\mathrm{d}x\right) \nonumber \\
 &\overset{\eqref{4}}{\leq}& cR^{-2}\left(\int_{B_{2R}-B_{R}}\Gamma_{1}^{\alpha+1-\frac{\mu_1}{2}}
 \,\mathrm{d}x+\int_{B_{2R}-B_{R}}\Gamma_{2}^{-\frac{\mu_2}{2}}\Gamma_{1}^{\alpha+1}\,\mathrm{d}x\right).
 \label{23}
\end{eqnarray}
Recall \eqref{5} and choose $\alpha$ according to
\begin{equation}
    \alpha\in\left(-1/2,\mathrm{min}\left\lbrace-1+\frac{\mu_1}{2},-1+\frac{\mu_2}{2}\right\rbrace\right).\label{26}
    \end{equation}
Here we note that -- depending on the values of $\mu_1$ and $\mu_2$ -- actually a negative exponent $\alpha$ can occur.
It follows from \eqref{9}  that
\begin{equation}
     cR^{-2}\left(\int_{B_{2R}-B_{R}}\Gamma_{1}^{\alpha+1-\frac{\mu_1}{2}}\,\mathrm{d}x+\int_{B_{2R}-B_{R}}
     \Gamma_{2}^{-\frac{\mu_2}{2}}\Gamma_{1}^{\alpha+1}\,\mathrm{d}x\right)\leq cR^{-2}\int_{B_{2R}-B_{R}}
     c\,\mathrm{d}x \leq c <\infty, \label{12}
\end{equation}
recalling that $c$ is independent of $R$.
Combining \eqref{12} and \eqref{23} it is obvious that (by passing to the limit $R\to\infty$)
\begin{equation}
    \int_{\mathbb{R}^2}D^2f(\nabla u)(\nabla\partial_1 u,\nabla\partial_1 u)\Gamma_{1}^{\alpha}\, \mathrm{d}x < \infty \label{13}
\end{equation}
for $\alpha$ satisfying \eqref{26}.\\

\newline
As in the proof of Proposition 6.1 from \cite{MR4502898} (with $l=1$) and by applying the Cauchy-Schwarz inequality we get
\begin{eqnarray}
 \lefteqn{\int_{B_{2R}}D^2f(\nabla u)(\nabla\partial_1 u,\nabla\partial_1 u)\eta^{2}\Gamma_{1}^{\alpha}\,\mathrm{d}x}\nonumber\\ 
 &   \leq  &\Biggl|\int_{B_{2R}}D^2 f(\nabla u)(\nabla\partial_{1}u,\nabla \eta^{2})\partial_{1}u\Gamma_{1}^{\alpha}
    \, \mathrm{d}x\Biggr| \nonumber \\ 
&\leq&c\left[\int_{\mathrm{spt}\nabla\eta}D^2f(\nabla u)(\nabla\partial_1 u,\nabla\partial_1 u)\eta^{2}\Gamma_{1}^{\alpha}
\,\mathrm{d}x\right]^{\frac{1}{2}} \left[\int_{\mathrm{spt}\nabla\eta}D^2f(\nabla u)(\nabla\eta,\nabla\eta)\Gamma_{1}^{\alpha+1}
\,\mathrm{d}x\right]^{\frac{1}{2}}.  \, \,\,\, \label{14}
\end{eqnarray}

The second integral on the right-hand side is bounded on account of our previous calculations.
Because of the validity of \eqref{13} the limit of the first integral for $R\to\infty$ is $0$. Thus \eqref{14} implies
\begin{equation}
    \nabla\partial_1 u\equiv 0. \label{15}
\end{equation}
In particular \eqref{15} guarantees the existence of a number $a\in\mathbb{R}$ such that
\begin{equation}
    \partial_1 u\equiv a. \label{16}
\end{equation}
From \eqref{16} we obtain
\begin{align*}
    u(x_1,x_2)-u(0,x_2)
    =\int_{0}^{x_1}\diff{}{t}u(t,x_2)\, \mathrm{dt}
    =\int_{0}^{x_1}a\,\mathrm{d}t
    =a\, x_1,
\end{align*}
implying
\begin{equation}
    u(x_1,x_2)=u(0,x_2)+a\, x_1. \label{17}
\end{equation}
Considering \eqref{16} again, equation \eqref{6} reduces to
\begin{equation}
    \partial_{2}\left(f_{2}'(\partial_{2}u)\right)=0. \label{18}
\end{equation}
We set $\varphi(t)\vcentcolon= u(0,t)$ and interpret the PDE \eqref{18} as the ODE
\begin{equation}
    \diff{}{t}\left(f_{2}'(\varphi'(t))\right)=0 \label{19}
\end{equation}
implying 
\begin{equation*}
f_{2}'(\varphi'(t)))=\mathrm{const.}
\end{equation*}
Since $f_{2}'$ is strictly monotonically increasing, this just means
\begin{equation*}
    \varphi '(t)=b , \quad t\in\mathbb{R},
\end{equation*}
for some real number $b$,
which consequently gives
\begin{equation}
    u(x_1,x_2)=a\, x_1 +b\, x_2 +c, \quad a,b,c \in\mathbb{R} \label{20}
\end{equation}
completing our proof. \qed

\section{Remaining proofs}

\textbf{ad Theorem 1.7.} Let the assumptions of Theorem 1.7 hold and assume w.l.o.g. that we have inequality \eqref{24} 
from Definition 1.6. Consider the mixed term in the last line of \eqref{23} and note that on account of \eqref{24} we may estimate
\begin{equation}
    \Gamma_{2}^{-\frac{\mu_2}{2}}\Gamma_1^{\alpha+1}\leq c\Gamma_{2}^{-\frac{\mu_2}{2}(1-2\rho(\alpha+1))}.
\end{equation}
The validity of $\rho<1$ allows us to choose $\alpha$ sufficiently close to $-1/2$ such that $1-2\rho(\alpha+1)>0$ 
which again yields \eqref{13} and allows us to proceed as before giving our claim. \qed \\

\textbf{ad Theorem 1.8.} Suppose by contradiction that there exists a real number $\hat{c}>0$ such that w.r.t the set $M_{1,\rho}$
\begin{equation}
    \hat{c}\leq\liminf_{|x_1|\to\infty}\frac{|\partial_1 u(x_1,x_2)|}{|x_1|^{\kappa_1}}. \label{31}
\end{equation}
For intervals $I_1,I_2\subset\mathbb{R}$ we let
\begin{equation*}
    S_{I_1;I_2}\vcentcolon=\lbrace x\in\mathbb{R}^2:|x_1|\in I_1, |x_2|\in I_2\rbrace.
\end{equation*}
We fix $0<R_1<R_2$ and consider 
\begin{equation*}
    \eta\in\mathit{C}_{0}^{\infty}(S_{[0,2R_1);[0,2R_2)}), \quad 0\leq\eta\leq 1, \quad \eta\equiv 1 \quad \text{on} 
    \quad S_{[0,2R_1);[0,2R_2)} \, ,
\end{equation*}
\begin{equation}
    \text{spt}\partial_1\eta\subset S_{(R_1,2R_1);[0,2R_2)} \, , \quad \text{spt}\partial_2\eta \subset S_{[0,2R_1);(R_2,2R_2)}\, , \label{32}
\end{equation}
\begin{equation}
    |\partial_1\eta|\leq c/R_1 \, , \quad |\partial_2\eta|\leq c/R_2. \label{33}
\end{equation}
Exactly as in \eqref{23} one obtains using \eqref{32} and \eqref{33}
\begin{eqnarray}
 \lefteqn{\int_{S_{[0,2R_1);[0,2R_2)}}D^2f(\nabla u)(\nabla \partial_1 u, \nabla \partial_1 u)\eta^2\Gamma_1^{\alpha} \, \mathrm{d}x} \nonumber\\ 
 &\leq& \int_{S_{[0,2R_1);[0,2R_2)}}D^2f(\nabla u)(\nabla \eta, \nabla \eta)\eta^2\Gamma_1^{1+\alpha} \, \mathrm{d}x \nonumber \\
  &  \leq &\frac{c}{R_1^2}\int_{S_{(R_1,2R_1);[0,2R_2)}}\Gamma_1^{\alpha+1-\frac{\mu_1}{2}} \,\mathrm{d}x
    +\frac{c}{R_2^2}\int_{S_{[0,2R_1);(R_2,2R_2)}}\Gamma_2^{-\frac{\mu_2}{2}}\Gamma_1^{\alpha+1} \, \mathrm{d}x. \label{34}
\end{eqnarray}
By definition we have 
\begin{equation*}
    |S_{[0,2R_1);(R_2,2R_2)}|\leq cR_1R_2.    
\end{equation*}
Moreover, our assumption $\kappa_2<1$ implies that $\alpha$ can be chosen such that in the case $\kappa_1>0$
\begin{equation}
    -\frac{1}{2}<\alpha<\frac{1}{2\kappa_2}-1. \label{35}
\end{equation}
In the case $\kappa_2=0$ we do not need an additional condition. We apply assumption \eqref{28}, which leads to
\begin{align}
    \frac{1}{R_2^2}\int_{S_{[0,2R_1);(R_2,2R_2)}}\Gamma_2^{-\frac{\mu_2}{2}}\Gamma_1^{\alpha+1} 
    \, \mathrm{d}x &\leq \frac{c}{R_2^2}\int_{S_{[0,2R_1);(R_2,2R_2)}}|x_2|^{2\kappa_2(\alpha+1)}\, \mathrm{d}x \nonumber \\
    &\leq cR_1R_2^{2\kappa_2(\alpha+1)-1}. \label{36}
\end{align}
Let us consider the first integral on the right-hand side of \eqref{34} recalling that $\alpha+1-\mu_1/2<0$. Assumption \eqref{31} implies 
\begin{align}
    \frac{1}{R_1^2}\int_{S_{(R_1,2R_1);[0,2R_2)}}\Gamma_1^{\alpha+1-\frac{\mu_1}{2}}\, \mathrm{d}x&\leq 
    \frac{c}{R_1^2}\int_{(R_1,2R_1);[0,2R_2)}|x_1|^{2\kappa_1(\alpha+1-\frac{\mu_1}{2})}\, \mathrm{d}x \nonumber \\
    &\leq cR_2R_1^{-\kappa_1(\mu_1-2(\alpha+1))-1}, \label{37}
\end{align}
where we suppose that $\mu_1>2(\alpha+1)$ by choosing $\alpha$ sufficiently close to $-1/2$. If we further suppose that
\begin{equation}
    R_2=R_1^{1+\rho} \quad \text{with a positive real number} \quad \rho<\kappa_1(\mu_1-1), \label{38}
\end{equation}
then by decreasing $\alpha$, if necessary, still satisfying $\alpha>-1/2$, we obtain from \eqref{37}
\begin{equation}
    \frac{1}{R_1^2}\int_{S_{(R_1,2R_1);[0,2R_2)}}\Gamma_1^{\alpha+1-\frac{\mu_1}{2}}\, \mathrm{d}x \to 0 
    \quad \text{as} \quad R_1\to\infty. \label{39}
\end{equation}
Using $R_2=R_1^{1+\rho}$ (recall \eqref{38}) we return to \eqref{36} recalling that 
by the choice \eqref{35} we have $2\kappa_2(\alpha+1)-1<0$. We calculate
\begin{equation}
    R_1R_1^{(1+\rho)(2\kappa_2(\alpha+1)-1)}=R_1^{2\kappa_2(\alpha+1)+\rho(2\kappa_2(\alpha+1)-1)}. \label{40}
\end{equation}
If we suppose that
\begin{equation}
    \kappa_2+\rho(\kappa_2-1)<0, \label{41}
\end{equation}
then we may choose $\alpha>-1/2$ sufficiently small such that the exponent on the right-hand side of \eqref{40} is negative, 
hence together with \eqref{36}
\begin{equation}
    \frac{1}{R_2^2}\int_{S_{[0,2R_1);(R_2,2R_2)}}\Gamma_2^{-\frac{\mu_2}{2}}\Gamma_1^{\alpha+1}\, \mathrm{d}x \to 0 
    \quad \text{as} \quad R_1\to \infty. \label{42}
\end{equation}
By \eqref{34}, \eqref{39} and \eqref{42} it follows that
\begin{equation}
    \int_{S_{[0,2R_1);[0,2R_2)}}D^2f(\nabla u)(\nabla\partial_1 u, \nabla\partial_1 u)\eta^2\Gamma_1^{\alpha}
    \, \mathrm{d}x \to 0 \quad \text{as} \quad R_1\to\infty, \label{43}
\end{equation}
provided that we have \eqref{38} and \eqref{41}, i.e. provided that we have \eqref{30} which is a consequence of \eqref{27}. 
Hence we have \eqref{43} which exactly as in the proof of Theorem 1.4 shows that $u$ has to be an affine function and 
this contradicts \eqref{31} which in turn proves Theorem 1.8. \qed\\

\textbf{ad Theorem 1.9.} W.l.o.g.~we suppose that $\Theta \in L^\infty(\mathbb{R}^2)$ and that $u \in \mathit{C}^3(\mathbb{R}^2)$,
$f_1$, $f_2 \in \mathit{C}^3(\mathbb{R})$. Otherwise we argue in a weak sense. Let
\[
w_i := f_i'(\partial_i u), \quad i=1,2 .
\]
Then we have
\begin{equation}\label{new 2a}
\partial_1 w_1 + \partial_2 w_2 = 0 \quad\mbox{on $\mathbb{R}^2$} ,
\end{equation}
hence
\begin{equation}\label{new 3}
\partial_{11}w_1 + \partial_1 \partial_2 w_2 =  \partial_{11}w_1 + \partial_2 \partial_1 w_2  = 0 \quad\mbox{on $\mathbb{R}^2$} .
\end{equation}
A direct calculation shows
\[
\partial_1 w_2 = \partial_1 \Big(f_2'(\partial_2 u)\Big) = f_2''(\partial_2 u) \partial_1 \partial_2 u = \Theta f_1''(\partial_1 u) \partial_2 \partial_1 u
= \Theta \partial_2 w_1 
\]
and the weak form of \eqref{new 3} reads as
\begin{equation}\label{new 4}
\int_{\mathbb{R}^2} \left(\begin{array}{c}
\partial_1 w_1\\ \Theta \partial_2 w_1
\end{array}\right) \cdot \nabla \varphi \dx = 0, \quad \varphi \in C^1_0(\mathbb{R}^2) .
\end{equation}
Inserting $\varphi = \eta^2 w_1$ with suitable $\eta \in C^1_0(B_{2R})$ such that $0 \leq \eta \leq 1$, $\eta \equiv 1$ on $B_R$
and $|\nabla \eta| \leq c/R$, we obtain
\begin{eqnarray}\label{new 5}
\lefteqn{\int_{B_{2R}} |\partial_1 w_1|^2 \eta^2 \dx + \int_{B_{2R}} \Theta |\partial_2 w_1|^2\eta^2 \dx}\nonumber\\
&= & - 2 \int_{B_{2R}-B_R} \eta \partial_1 w_1\, \partial_1 \eta \, w_1 \dx 
- 2 \int_{B_{2R}-B_R} \Theta\, \eta\,  \partial_2 w_1\, \partial_2 \eta \,w_1 \dx .
\end{eqnarray}
Applying Young's inequality and using $w_1$, $\Theta \in L^{\infty}(\mathbb{R}^2)$ we obtain that
\begin{equation}\label{new 6}
\int_{\mathbb{R}^2} \Big( |\partial_1 w_1|^2 + \Theta |\partial_2 w_1|^2\Big) \dx < \infty .
\end{equation}
We then return to \eqref{new 5} and apply the inequality of Cauchy-Schwarz to obtain
\begin{eqnarray}\label{new 7}
\lefteqn{\int_{B_{2R}} |\partial_1 w_1|^2 \eta^2 \dx + \int_{B_{2R}} \Theta |\partial_2 w_1|^2 \eta^2\dx}\nonumber\\
&\leq & 2 \Bigg[\int_{B_{2R}-B_R}\eta^2 |\partial_1 w_1|^2 \dx\Bigg]^{\frac{1}{2}} 
\Bigg[\int_{B_{2R}-B_R} |\partial_1 \eta|^2 w_1^2 \dx \Bigg]^{\frac{1}{2}}\nonumber\\
&&+ 2 \Bigg[ \int_{B_{2R}-B_R} \Theta \eta^2 |\partial_2 w_1|^2 \dx\Bigg]^{\frac{1}{2}}
\Bigg[ \int_{B_{2R}-B_R} \Theta |\partial_2 \eta|^2 w_1^2 \dx \Bigg]^{\frac{1}{2}} .
\end{eqnarray}
On the right-hand side of \eqref{new 7} we observe that for both parts the first integral is vanishing when passing to the limit $R\to \infty$
since we have \eqref{new 6}, while the remaining integrals stay uniformly bounded.\\

This gives $\partial_1 w_1 = 0$ and $\partial_2 w_1 = 0$ since we have $\Theta >0$. Hence we obtain $w_1 \equiv c_1$
for some constant $c_1$. The monotonicity of $f_1'$ then implies $\partial_1 u \equiv  \tilde{c}_1$
for some different constant $\tilde{c}_1$.\\

By \eqref{new 2a} we then also have $\partial_2 w_2 =0$. Since we have already observed above that $\partial_1 w_2 = \Theta \partial_2 w_1$, we
deduce $\partial_2 u \equiv \tilde{c}_2$ for some other real number $\tilde{c}_2$ and 
in conclusion $u$ must be an affine function which completes the proof of
Theorem 1.9. \qed\\

\bibliographystyle{unsrtnat}
\bibliography{SmallRemarkBernstein} 

\begin{thebibliography}{9}
\providecommand{\natexlab}[1]{#1}
\providecommand{\url}[1]{\texttt{#1}}
\expandafter\ifx\csname urlstyle\endcsname\relax
  \providecommand{\doi}[1]{doi: #1}\else
  \providecommand{\doi}{doi: \begingroup \urlstyle{rm}\Url}\fi

\bibitem[Bernstein(1927)]{MR1544873}
Serge Bernstein.
\newblock \"{U}ber ein geometrisches {T}heorem und seine {A}nwendung auf die
  partiellen {D}ifferentialgleichungen vom elliptischen {T}ypus.
\newblock \emph{Math. Z.}, 26\penalty0 (1):\penalty0 551--558, 1927.

\bibitem[Dierkes et~al.(2010{\natexlab{a}})Dierkes, Hildebrandt, and
  Tromba]{MR2778928}
Ulrich Dierkes, Stefan Hildebrandt, and Anthony~J. Tromba.
\newblock \emph{Global analysis of minimal surfaces}, volume 341 of
  \emph{Grundlehren der mathematischen Wissenschaften [Fundamental Principles
  of Mathematical Sciences]}.
\newblock Springer, Heidelberg, second edition, 2010{\natexlab{a}}.

\bibitem[Dierkes et~al.(2010{\natexlab{b}})Dierkes, Hildebrandt, and
  Tromba]{MR2760441}
Ulrich Dierkes, Stefan Hildebrandt, and Anthony~J. Tromba.
\newblock \emph{Regularity of minimal surfaces}, volume 340 of
  \emph{Grundlehren der mathematischen Wissenschaften [Fundamental Principles
  of Mathematical Sciences]}.
\newblock Springer, Heidelberg, second edition, 2010{\natexlab{b}}.
\newblock With assistance and contributions by A. K\"{u}ster.

\bibitem[Dierkes et~al.(2010{\natexlab{c}})Dierkes, Hildebrandt, and
  Sauvigny]{MR2566897}
Ulrich Dierkes, Stefan Hildebrandt, and Friedrich Sauvigny.
\newblock \emph{Minimal surfaces}, volume 339 of \emph{Grundlehren der
  mathematischen Wissenschaften [Fundamental Principles of Mathematical
  Sciences]}.
\newblock Springer, Heidelberg, second edition, 2010{\natexlab{c}}.
\newblock With assistance and contributions by A. K\"{u}ster and R. Jakob.

\bibitem[Nitsche(1989)]{MR1015936}
Johannes C.~C. Nitsche.
\newblock \emph{Lectures on minimal surfaces. {V}ol. 1}.
\newblock Cambridge University Press, Cambridge, 1989.
\newblock Introduction, fundamentals, geometry and basic boundary value
  problems, Translated from the German by Jerry M. Feinberg,.

\bibitem[Nitsche(1957)]{MR90833}
Johannes C.~C. Nitsche.
\newblock Elementary proof of {B}ernstein's theorem on minimal surfaces.
\newblock \emph{Ann. of Math. (2)}, 66:\penalty0 543--544, 1957.

\bibitem[Osserman(1986)]{MR852409}
Robert Osserman.
\newblock \emph{A survey of minimal surfaces}.
\newblock Dover Publications, Inc., New York, second edition, 1986.

\bibitem[Bildhauer and Fuchs(2020)]{MR4151290}
Michael Bildhauer and Martin Fuchs.
\newblock Splitting type variational problems with linear growth conditions.
\newblock \emph{J. Math. Sci. (N.Y.) 250(2), 2020. Problems in mathematical
  analysis}, 250\penalty0 (105):\penalty0 232--249, 2020.

\bibitem[Bildhauer and Fuchs(2022)]{MR4502898}
Michael Bildhauer and Martin Fuchs.
\newblock On the global regularity for minimizers of variational integrals:
  splitting-type problems in 2{D} and extensions to the general anisotropic
  setting.
\newblock \emph{J. Elliptic Parabol. Equ.}, 8\penalty0 (2):\penalty0 853--884,
  2022.

\end{thebibliography}

\end{document}